\documentclass{birkjour}
\usepackage{graphicx}
\usepackage{amsfonts, amssymb, amsmath, mathrsfs, latexsym}
\usepackage[all]{xy}
\usepackage[dvipsnames]{xcolor}
\usepackage{mathtools}
\usepackage{tikz}
\usetikzlibrary{positioning}
\usepackage{xcolor}
\usepackage[normalem]{ulem}

\ifpdf
	\usepackage[unicode,pdftex]{hyperref}
	\input glyphtounicode
\else
	\usepackage[unicode,dvipdfm]{hyperref}
\fi

\newtheorem{thm}{Theorem}[section]
\newtheorem{prop}[thm]{Proposition}
\newtheorem{lem}[thm]{Lemma}

\theoremstyle{definition}

\newtheorem{rmk}[thm]{Remark}
\newtheorem*{mcs}{Modified Assertion of Severi}

\newcommand{\PP}{\ensuremath{\mathbb{P}}}

\begin{document}

\title[On  the Hilbert scheme of  linearly normal curves]
{On  the Hilbert scheme of  linearly normal curves in $\mathbb{P}^4$ of degree $d = g+1$  and genus $g$}

\thanks{This paper was prepared for publication when the first named author was enjoying hospitality of the Institute of Mathematics - Academia Sinica (Taiwan) to which he is grateful for the support and the stimulating atmosphere.   Both authors were supported in part by National Research Foundation of South Korea 
 (2017R1D1A1B031763).} 

\dedicatory{In memory of Professor R.D.M. Accola}

\author[Changho Keem]{Changho Keem}
\address{
Department of Mathematics,
Seoul National University\\
Seoul 151-742,  
South Korea}

\email{ckeem1@gmail.com}

\thanks{}

\author[Yun-Hwan Kim]{Yun-Hwan Kim}
\address{Department of Mathematics,
Seoul National University\\
Seoul 151-742, 
South Korea}
\email{yunttang@snu.ac.kr}

\subjclass{Primary 14C05, Secondary 14H10}

\keywords{Hilbert scheme, algebraic curves, linear series}

\date{\today}
\maketitle

\begin{abstract}
We denote by $\mathcal{H}_{d,g,r}$ the Hilbert scheme of smooth curves, which is the union of components whose general point corresponds to a smooth irreducible and non-degenerate curve of degree $d$ and genus $g$ in $\PP^r$. 
In this article, we show that any non-empty $\mathcal{H}_{g+1,g,4}$ has only one component whose general element is linearly normal  unless $g=9$.  If $g=9$, we show that $\mathcal{H}_{g+1,g,4}$ is reducible with two components and a general element of each component is linearly normal.
This establishes the validity of a certain modified version of an assertion of Severi regarding the irreducibility of $\mathcal{H}_{d,g,r}$ for the case $d=g+1$ and $r=4$.

\end{abstract}
\section{\quad An overview, preliminaries and basic set-up}

Given non-negative integers $d$, $g$ and $r\ge 3$, let $\mathcal{H}_{d,g,r}$ be the Hilbert scheme of smooth curves parametrizing smooth irreducible and non-degenerate curves of degree $d$ and genus $g$ in $\PP^r$.

After Severi asserted with an incomplete proof that $\mathcal{H}_{d,g,r}$ is irreducible for $d\ge g+r$  in \cite{Sev}, the irreducibility of $\mathcal{H}_{d,g,r}$ has been studied by several authors. Ein proved Severi's claim for $r=3$ \& $r=4$; cf.  \cite[Theorem 4]{E1} and \cite[Theorem 7]{E2}.

For families of curves of lower degrees in $\mathbb{P}^3$, there are several works due to many people. The most updated result is that any non-empty $\mathcal{H}_{d,g,3}$ is irreducible for every $d\ge g$; cf. \cite[Theorem 1.5]{KK},  \cite[Proposition 2.1 and Proposition 3.2 ]{KKL}, \cite[Theorem 3.1]{I} and \cite{KKy1}.

For families of curves in $\PP^4$ of lower degree $d\le g+3$, Hristo Iliev proved the irreducibility of $\mathcal{H}_{d,g,4}$ for $d=g+3$, $g\ge 5$ and $d=g+2$, $g\ge 11$; cf. \cite{I}. Quite recently,  there has been a  minor extension of the result of Hristo Iliev regarding the irreducibility of $\mathcal{H}_{g+2,g,4}$ for low genus cases; i.e.  $\mathcal{H}_{g+2,g,4}$ is irreducible  and generically reduced for any genus $g$ as long as $\mathcal{H}_{g+2,g,4}$ is non-empty; cf. \cite[Corollary 2.2]{KKy2}. In this article we study the next case $\mathcal{H}_{g+1,g,4}$ - the Hilbert scheme of curves in $\mathbb{P}^4$ of degree $g+1$ - which has not been studied or understood well enough before.

On the other hand,  it is quite well known that Severi's assertion turned out to be untrue for curves in higher dimensional projective space $\mathbb{P}^r$ with $r\ge 5$; cf. \cite[Theorem 2.3]{Keem} for the irreducibility of $\mathcal{H}_{2g-8, g, r}$ for $r$ in the range $\frac{2g-7}{3} \le r \le g-8$ and \cite{CS} for several other variations of it.

Before proceeding, it should be remarked that there are several sources  in the literature suggesting that what Severi indeed had in mind in his original assertion regarding the irreducibility of $\mathcal{H}_{d,g,r}$ might have been the following statement, which we call the {\bf  Modified 
Assertion of Severi} and this is one of the motivation of our study in this paper;  e.g. cf. \cite[page 489]{CS} or AMS Mathematical Reviews MR1221726.

\begin{mcs}  A nonempty $\mathcal{H}^\mathcal{L}_{d,g,r}$ is irreducible for any triple $(d,g,r)$ in the Brill-Noether range  $\rho (d,g,r):=g-(r+1)(g-d+r)\ge 0$ , where $\mathcal{H}^\mathcal{L}_{d,g,r}$ is the union of those components of ${\mathcal{H}}_{d,g,r}$ whose general element is linearly normal.
\end{mcs}
\begin{rmk}\label{review}

(1)
As far as the authors know, all the examples of Hilbert schemes of smooth curves violating the Severi's original assertion in the Brill-Noether range  are those with at least one extra component whose general element is a non-linearly normal curve, i.e. a curve embedded in the corresponding projective space by an incomplete linear system.
In other words, in all such known examples it holds that $\mathcal{H}^\mathcal{L}_{d,g_r}\neq\mathcal{H}_{d,g,r}$ whereas   $\mathcal{H}^\mathcal{L}_{d,g,r}$ is still irreducible. 

\vskip 4pt
\noindent
(2)
By the irreducibility of $\mathcal{H}_{d,g,3}$ in the range $d\ge g$ and the proofs therein, $\mathcal{H}^\mathcal{L}_{d,g,3}=\mathcal{H}_{d,g,3}$ is irreducible either inside or outside the Brill-Noether range as long as $\mathcal{H}^\mathcal{L}_{d,g,3}$ is non-empty.
Likewise,  by the irreducibility of $\mathcal{H}_{d,g,4}$ in the range $d\ge g+2$, $\mathcal{H}^\mathcal{L}_{d,g,4}=\mathcal{H}_{d,g,4}$ is irreducible.

\vskip 4pt
\noindent
(3)
There is a  reducible example of $\mathcal{H}_{d,g,4}$ such that $\mathcal{H}^\mathcal{L}_{d,g,4}=\mathcal{H}_{d,g,4}$ outside the Brill-Noether range, which incidentally occurs in the case of our study $d=g+1$ and $g=9$. This is one of rare examples of reducible 
 $\mathcal{H}_{d,g,r}^{\mathcal{L}}$ when $r=4$.; Theorem \ref{lowgenus} (2).

\vskip 4pt
\noindent
(4)
For the Hilbert scheme of curves in $\mathbb{P}^3$, there is  also an example of reducible Hilbert scheme $\mathcal{H}_{d,g,3}$ such that $\mathcal{H}^\mathcal{L}_{d,g,3}=\mathcal{H}_{d,g,3}$ outside the Brill-Noether range which  occurs in case $d=g-1$. In fact it is easy to prove that  $\mathcal{H}_{9,10,3}=\mathcal{H}_{9,10,3}^\mathcal{L}$ is reducible with exactly two components, one consisting of trigonal curves lying on a smooth quadric and the other one consisting of complete intersection of two cubics both of the expected dimension 36; cf. Remark \ref{genus9}. Note for all $d\ge g$, the  Hilbert scheme $\mathcal{H}_{d,g,3}$ is irreducible.

\vskip 4pt
\noindent
(5)
There is a reducible example $\mathcal{H}_{d,g,4}$  of another kind; $\mathcal{H}^\mathcal{L}_{d,g,4}\neq\mathcal{H}_{d,g,4}$ but {\bf $\mathcal{H}^\mathcal{L}_{d,g,4}$} itself is irreducible. This also occurs in the case of our study  $d=g+1$ and $g=12$, which is outside the Brill-Noether range; Theorem \ref{lowgenus} (4).

\vskip 4pt
\noindent
(6)
In view of all these, the Modified Assertion of Severi clearly makes sense at least for $r=4$.  In particular the irreducibility of $\mathcal{H}^\mathcal{L}_{g+1,g,4}$  in the range $\rho (g+1,g,4)=g-15\ge 0$ seems to be worthwhile to be studied along this line.

\vskip 4pt
\noindent
(7)
We call the locus $\mathcal{H}^\mathcal{L}_{d,g,r}$ the {\bf Hilbert scheme of linearly normal curves}.

\vskip 4pt
\noindent
(8)
The irreducibility of $\mathcal{H}^\mathcal{L}_{d,g,r}$ in the Brill-Noether range $\rho (d,g,r)\ge 0$ was also suggested by J. Harris (implicitly or indirectly in a sense) in \cite[ page 77]{H1}, where he conjectured that ${\mathcal{H}}_{d,g,r}$ is irreducible in the Brill-Noether range. Of course,  the conjecture he made there - which is even stronger than the original assertion of Severi - is not valid by those examples violating Severi's assertion. However it also looks like that J. Harris might have meant the irreducibility of the Hilbert scheme of linearly normal curves 
 $\mathcal{H}^\mathcal{L}_{d,g,r}$, not the whole $\mathcal{H}_{d,g,r}$.

\vskip 4pt
\noindent
(9)
Needless to say, in case $d\ge g+r+1$ there is no  complete linear system of degree $d$ of dimension $r$ by the Riemann-Roch formula. Therefore in this range we have $\mathcal{H}^\mathcal{L}_{d,g,r}=\emptyset$ and hence the Modified Assertion of Severi makes sense only if $g-d+r\ge 0$.

\vskip 4pt
\noindent
(10)
For any $(d,g,r)$ with $r\ge 3$ in the Brill-Noether range, the author does not know of any example of reducible Hilbert scheme of linearly normal curves 
$\mathcal{H}^\mathcal{L}_{d,g,r}$.
Furthermore the authors have a strong feeling that the Modified Assertion of Severi would hold for any $r\ge 3$ at least inside the Brill-Noether range with $g-d+r\ge 0$. As the first attempt toward such an extensive settlement,
we deal with the particular case $d=g+1$ and  $r=4$, which we believe is worthy of being studied 

\end{rmk}

The organization of this paper is as follows. After we briefly mention and recall several basic preliminaries in the remainder of this section, we start the next section with the two examples of $\mathcal{H}^\mathcal{L}_{g+1,g,4}$ for $g=9$ and $g=12$ which we mentioned before in the Remark \ref{review} (3) \& (5). As we shall see shortly, $\mathcal{H}^\mathcal{L}_{g+1,g,4}$ is reducible with two components of different dimensions for $g=9$. However  $\mathcal{H}^\mathcal{L}_{g+1,g,4}$ is irreducible when $g=12$ even though $\mathcal{H}_{g+1,g,4}$ itself is reducible with two components of the same dimension. We also deal with the irreducibility of $\mathcal{H}^\mathcal{L}_{g+1,g,4}$ for some low genus $g$, e.g.  $g=10$ or $g=11$. 

In the last part of the next section we finish the proof of our main
result by using Lemma \ref{dualbirva} which characterizes the residual series with respect to the canonical series of the complete linear series corresponding to the linearly normal curves under consideration. We use  the irreducibility of the Severi variety of plane curves in the final stage of the proof; cf. \cite{AC2} and~\cite{H2}.

For notations and conventions, we usually follow those in \cite{ACGH} and \cite{ACGH2}; e.g. $\pi (d,r)$ is the maximal possible arithmetic genus of an irreducible and non-degenerate curve of degree $d$ in $\PP^r$. Throughout we work over the field of complex numbers.

Before proceeding, we recall several related results which are rather well-known; cf. \cite{ACGH2}.
Let $\mathcal{M}_g$ be the moduli space of smooth curves of genus $g$. For any given isomorphism class $[C] \in \mathcal{M}_g$ corresponding to a smooth irreducible curve $C$, there exist a neighborhood $U\subset \mathcal{M}_g$ of the class $[C]$ and a smooth connected variety $\mathcal{M}$ which is a finite ramified covering $h:\mathcal{M} \to U$, as well as  varieties $\mathcal{C}$, $\mathcal{W}^r_d$ and $\mathcal{G}^r_d$ proper over $\mathcal{M}$ with the following properties:
\begin{enumerate}
\item[(1)] $\xi:\mathcal{C}\to\mathcal{M}$ is a universal curve, i.e. for every $p\in \mathcal{M}$, $\xi^{-1}(p)$ is a smooth curve of genus $g$ whose isomorphism class is $h(p)$,
\item[(2)] $\mathcal{W}^r_d$ parametrizes the pairs $(p,L)$ where $L$ is a line bundle of degree $d$ and $h^0(L) \ge r+1$ on $\xi^{-1}(p)$,
\item[(3)] $\mathcal{G}^r_d$ parametrizes the couples $(p, \mathcal{D})$, where $\mathcal{D}$ is possibly an incomplete linear series of degree $d$ and dimension $r$ on $\xi^{-1}(p)$ - which is usually denoted by $g^r_d$. 
\end{enumerate}

\vskip 4pt
Let $\widetilde{\mathcal{G}}$ ($\widetilde{\mathcal{G}}_\mathcal{L}$ resp.) be  the union of components of $\mathcal{G}^{r}_{d}$ whose general element $(p,\mathcal{D})$ of $\widetilde{\mathcal{G}}$ ($\widetilde{\mathcal{G}}_\mathcal{L}$ resp.) corresponds to a very ample (very ample and complete resp.) linear series $\mathcal{D}$ on the curve $C=\xi^{-1}(p)$. Note that an open subset of $\mathcal{H}_{d,g,r}$ consisting of points corresponding to smooth irreducible and non-degenerate curves is a $\mathbb{P}\textrm{GL}(r+1)$-bundle over an open subset of $\widetilde{\mathcal{G}}$. Hence the irreducibility of $\widetilde{\mathcal{G}}$ guarantees the irreducibility of $\mathcal{H}_{d,g,r}$. Likewise, the irreducibility of $\widetilde{\mathcal{G}}_\mathcal{L}$ ensures the irreducibility of 
 $\mathcal{H}_{d,g,r}^\mathcal{L}$.

We also make a note of the following well-known facts regarding the schemes $\mathcal{G}^{r}_{d}$ and $\mathcal{W}^r_d$; cf. \cite[Proposition 2.7, 2.8]{AC2}, \cite[2.a]{H1} and \cite[Theorem 1]{EH}.

\begin{prop}\label{facts}
For non-negative integers $d$, $g$ and $r$, let $\rho(d,g,r):=g-(r+1)(g-d+r)$ be the Brill-Noether number.
	\begin{enumerate}
	\item[\rm{(1)}] The dimension of any component of $\mathcal{G}^{r}_{d}$ is at least $3g-3+\rho(d,g,r)$ which is denoted by $\lambda(d,g,r)$. Moreover, if $\rho(d,g,r)\geq0$, there exists a unique component $\mathcal{G}_0$ of $\widetilde{\mathcal{G}}$ which dominates $\mathcal{M}$(or $\mathcal{M}_g$).

	\item[\rm{(2)}] Suppose $g>0$ and let $X$ be a component of $\mathcal{G}^{2}_{d}$ whose general element $(p,\mathcal{D})$ is such that $\mathcal{D}$ is a birationally very ample linear series on $\xi^{-1}(p)$. Then 
	\[\dim X=3g-3+\rho(d,g,2)=3d+g-9.\]
		\end{enumerate}
\end{prop}
\begin{rmk}\label{principal}

(1)
In the Brill-Noether range, the unique component $\mathcal{G}_0$ of $\widetilde{\mathcal{G}}$ (and the corresponding component $\mathcal{H}_0$ of $\mathcal{H}_{d,g,r}$ as well) which dominates $\mathcal{M}$ or $\mathcal{M}_g$ is called the  ``principal component".  

\vskip 4pt
\noindent
(2) In the range $d\le g+r$ inside the Brill-Noether range, the principal component $\mathcal{G}_0$ which has the expected dimension is one of the components of $\widetilde{\mathcal{G}}_\mathcal{L}$ (cf. \cite[2.1 page 70]{H1}),   and therefore $\widetilde{\mathcal{G}}_\mathcal{L}$ or 
 $\mathcal{H}_{d,g,r}^\mathcal{L}$ is non-empty. Hence what we are chasing after is if $\mathcal{G}_0$ is the only component of 
 $\widetilde{\mathcal{G}}_\mathcal{L}$.
\end{rmk}

We recall that the family of plane curves of degree $d$ in $\mathbb{P}^2$ are naturally parametrised by the projective space $\mathbb{P}^N$, $N=\frac{d(d+3)}{2}$. Let $\Sigma_{d,g}\subset \PP^N$ be the Severi variety of plane curves of degree $d$ with geometric genus $g$. We also recall that a general point of $\Sigma_{d,g}$ corresponds to an  irreducible plane curve of degree $d$ having $\delta:=\frac{(d-1)(d-2)}{2}-g$ nodes and no other singularities. 
The following theorem of Harris is fundamental; cf. \cite[Theorem 10.7 and 10.12]{ACGH2} or \cite[Lemma 1.1, 1.3 and 2.3]{H2} .

\begin{thm}\label{severi}
$\Sigma_{d,g}$ is irreducible of dimension $3d+g-1=\lambda(d,g,2)+\dim\mathbb{P}\rm{GL}(3)$. 
\end{thm}

Denoting by $\mathcal{G'}\subset \mathcal{G}^{2}_{d}$ the 
union of components whose general element $(p,\mathcal{D})$  is such that $\mathcal{D}$ is birationally very ample  on $C=\xi^{-1}(p)$, we 
remark that an open subset of the Severi variety $\Sigma_{d,g}$ is a $\mathbb{P}\textrm{GL}(3)$-bundle over an open subset of  $\mathcal{G}'$.  Therefore, as  an immediate consequence of Theorem \ref{severi}, the irreducibility of  $\Sigma_{d,g}$ implies the 
irreducibility of the locus $\mathcal{G'}\subset \mathcal{G}^{2}_{d}$ and vice versa. We make a note of this  observation as the following lemma.

\begin{lem}\label{Gisirred}
Let $\mathcal{G}'\subset \mathcal{G}^{2}_{d}$ be the union of components whose general element $(p,\mathcal{D})$ is such that $\mathcal{D}$ is birationally very ample on $C=\xi^{-1}(p)$. Then $\mathcal{G}'$ is irreducible.
\end{lem}

We will utilize the  following upper bound of the dimension of an irreducible component of $\mathcal{W}^r_d$, which was proved  and used effectively in \cite{I}.

\begin{prop}[\rm{\cite[Proposition 2.1]{I}}]\label{wrdbd}
Let $d,g$ and $r\ge 2$ be positive integers such that  $d\le g+r-2$ and let $\mathcal{W}$ be an irreducible component of $\mathcal{W}^{r}_{d}$. For a general element $(p,L)\in \mathcal{W}$, let $b$ be the degree of the base locus of the line bundle $L=|D|$ on $C=\xi^{-1}(p)$. Assume further that for a general $(p,L)\in \mathcal{W}$ the curve $C=\xi^{-1}(p)$ is not hyperelliptic. If the moving part of $L=|D|$ is
	\begin{itemize}
	\item[\rm{(a)}] very ample and $r\ge3$, then 
	$\dim \mathcal{W}\le 3d+g+1-5r-2b$;
	\item[\rm{(b)}] birationally very ample, then 
	$\dim \mathcal{W}\le 3d+g-1-4r-2b$;
	\item[\rm{(c)}] compounded, then 
	$\dim \mathcal{W}\le 2g-1+d-2r$.
	\end{itemize}
\end{prop}

\section{\quad Irreducibility of $\mathcal{H}^\mathcal{L}_{g+1,g,4}$}
The main result of this article is the following theorem, from which the Modified Severi's Assertion (for $d=g+1$ and $r=4$) follows immediately. The result is stronger than what we have expected.
 
\begin{thm}\label{H_g+1,g,4}
Every non-empty $\mathcal{H}_{g+1,g,4}^\mathcal{L}$ is irreducible unless $g=9$.
\end{thm}

We begin with making a note of the following facts when the genus of curves under consideration is relatively low.

\begin{thm}\label{lowgenus}
	\begin{enumerate}
		\item[\rm{(1)}] $\mathcal{H}_{g+1,g,4}=\mathcal{H}_{g+1,g,4}^\mathcal{L}=\emptyset$ for $g\leq8$.		
		\item[\rm{(2)}] For $g=9$, $\mathcal{H}_{10,9,4}=\mathcal{H}_{10,9,4}^\mathcal{L}$ is reducible with two components of dimensions $42$ and $43$.
		\item[\rm{(3)}] For $g=10$,  $\mathcal{H}_{11,10,4}=\mathcal{H}_{11,10,4}^\mathcal{L}$ is irreducible of the expected dimension $46$.
		\item[\rm{(4)}] For $g=12$, $\mathcal{H}_{13,12,4}$ is reducible with two components of the same expected dimension $54$, whereas 
		$\mathcal{H}_{13,12,4}^\mathcal{L}$ is irreducible. 
	\end{enumerate}
\end{thm}
\noindent
{\it Proof.}
(1) For $0\le g\le 2$, every curve of degree $g+1$ in $\mathbb{P}^4$ is degenerate. For $3\le g \le 8$, one easily computes that  $g\le \pi (g+1,4) < g$ by the Castelnuovo genus bound, a contradiction. 

\vskip 4pt
\noindent
(2) A smooth curve of genus $g=9$ in $\mathbb{P}^4$ of degree $10$ is an extremal curve lying on a cubic scroll $S$; 
$\pi (10,4)=9$.   We remark that the reducibility of  $\mathcal{H}_{10,9,4}$ is a very special case of the result by C. Ciliberto regarding the Hilbert scheme of curves of maximal genus; cf. \cite[Theorem 1.4(ii)]{CC}. However we carry out an elementary dimension count for the convenience of readers. Let $H$ and $L$ be classes of a hyperplane and  a ruling of the scroll $S$ and let $C\in|\alpha H+\beta L|$.  By computing the arithmetic genus of the $C$ lying on the scroll $S$, we see that either (i) $C\in |3H+L|$ or (ii) $C\in |4H-2L|$. 

\vskip 4pt
\noindent
(2-i) Let $\mathcal{H}_1$ be the component whose general member  $C\in\mathcal{H}_1$ is in the linear system $|3H+L|$ on a  scroll $S\subset\mathbb{P}^4$.
Since $(3H+L)\cdot L=3$, the rulings $|L|$ cut out a unique trigonal pencil $g^1_3$ on $C$.  Furthermore, the residual series $|K_C-\mathcal{D}|$ of the hyperplane series $|\mathcal{D}|$ is cut out on $C$ by the series $|K_S+C-H|=|(-2H+L)+(3H+L)-H|=|2L|$ and hence $|K_C-\mathcal{D}|=2g^1_3$. Conversely any trigonal curve of genus $9$ has a (unique) very ample $g^4_{10}$ which is of the form $|K_C-2g^1_3|$. Therefore it follows  that $\mathcal{H}_1$ is a $\mathbb{P}GL(5)$-bundle over an irreducible locus $\mathcal{G}_1\subset\mathcal{G}^4_{10}$ which consists of the residual series of $2g^1_3$ on  trigonal curves. 
Hence we have $$\dim\mathcal{H}_1=\dim\mathcal{G}_1+\dim\mathbb{P}GL(5)=\dim\mathcal{M}^1_{g,3}+\dim\mathbb{P}GL(5)=43.$$  
\noindent
(2-ii) Let $\mathcal{H}_2$ be the component whose general member  $C\in\mathcal{H}_2$
is a curve in the linear system  $|4H-2L|$  on a scroll $S\subset\mathbb{P}^4$. We note that $C$ is $4$-gonal and is not trigonal by the Castelnuovo-Severi inequality; cf.  \cite[Theorem 3.5]{Accola}.
We also note  that the residual series $|K_C-\mathcal{D}|=g^2_6$ of the hyperplane series $\mathcal{D}$ is cut out on $C$ by the series $|K_S+C-H|=|(-2H+L)+(4H-2L)-H|=|H-L|$ and hence $|K_C-\mathcal{D}|$ is  base-point-free and birationally very ample. It then follows that  $C$ is birational to a plane sextic with a one node or a cusp. Conversely a curve of genus $9$ which is birational to a plane sextic with a nodal singularity can be embedded in $\mathbb{P}^4$ as a curve of degree $10$ by the linear series $|g^2_6+g^1_4|$, where $g^1_4$ is the unique pencil of degree $4$ cut out by lines through the nodal singularity. Therefore we may conclude that $\mathcal{H}_2$ is a $\mathbb{P}GL(5)$-bundle over an irreducible component $\mathcal{G}_2\subset\mathcal{G}^4_{10}$ which consists of the residual series of  birationally very ample $g^2_6$'s. Hence $$\dim\mathcal{H}_2=\dim\mathcal{G}_2+\dim\mathbb{P}GL(5)=3\cdot 6+g-9+24=42$$ by Proposition \ref{facts}(2). We also see that  $\mathcal{H}_1\neq \mathcal{H}_2$ by  semi-continuity. A general element of  $\mathcal{H}_1$ or  $\mathcal{H}_2$ is an extremal curve in $\mathbb{P}^4$ and  hence corresponds to a  linearly normal curve so that  $\mathcal{H}^\mathcal{L}_{10,9,4}=\mathcal{H}_{10,9,4}=\mathcal{H}_1\cup\mathcal{H}_2$.

\vskip 4pt
\noindent
(3)
For $g=10$, we let ${\mathcal{G}}\subset \mathcal{G}^{4}_{g+1}$ be an irreducible component whose general element $(p, \mathcal{D})$ is a very ample linear series $\mathcal{D}$ on the curve $C=\xi^{-1}(p)$. We note that a smooth curve $C$ of genus $g=10$ cannot have an incomplete very ample $\mathcal{D}=g^4_{11}$ by Castelnuovo genus bound. Hence the residual series $|K_C-\mathcal{D}|$ is a complete $g^2_7$. 
If $|K_C-\mathcal{D}|=g^2_7$ has a base point, then $C$ is either birational to a smooth plane sextic or trigonal. In the former case, by counting the family whose residual series is of the form $|g^2_6+q|$, where $g^2_6$ is very ample, we have 
\begin{align*}\dim\mathcal{G}&=\dim\mathbb{P}(H^0(\mathbb{P}^2, \mathcal{O}(6)))-\dim\mathbb{P}\rm{GL}(3)+1=20\\
&<\lambda(11,10,4)=22
\end{align*}
a contradiction; cf. Proposition \ref{facts} (1). If $C$ is trigonal, the residual series of $|2g^1_3+q|=g^2_7$ is not even very ample; note that $|\mathcal{D}-r-s|=|K_C-2g^1_3-q-r-s|=
|K_C-3g^1_3|=g^3_9$ where $q+r+s\in g^1_3$.
Therefore it follows that  $|K_C-\mathcal{D}|=g^2_7$ is base-point-free and hence $C$ is birational to a (singular) plane septic. In this case we see that the component $\mathcal{G}$ is birational to the irreducible $\mathcal{G}'\subset\mathcal{G}^2_7$, which is the locus over which the Severi variety $\Sigma_{7,10}$ lies; Theorem \ref{severi} and  Lemma \ref{Gisirred}. Note that 
$\mathcal{H}_{11,10,4}\neq\emptyset$ since there exists a smooth curve in $\mathbb{P}^4$ of genus $10$ and degree $11$; e.g.  a proper transformation of a plane septic with $5$ nodal singularities on a del-Pezzo surface in $\mathbb{P}^4$,  blown up at five points in general position in $\mathbb{P}^2$; cf. \cite[Theorem 1.0.1]{rathman}. Hence we see that $\mathcal{H}_{11,10,4}=\mathcal{H}^\mathcal{L}_{11,10,4}$ is irreducible and 
\begin{align*}
\dim\mathcal{H}_{11,10,4}&=\dim\mathcal{G}+\dim\mathbb{P}GL(5)=\dim\mathcal{G}'+\dim\mathbb{P}GL(5)\\&=3g-3+\rho (7,10,2)+\dim\mathbb{P}GL(5)\\
&=\lambda(7,10,2)+\dim\mathbb{P}GL(5)=46.
\end{align*}

\vskip 4pt
\noindent
(4)
For $g=12$,  we also let ${\mathcal{G}}\subset \mathcal{G}^{4}_{g+1}$ be an irreducible component whose general element $(p, \mathcal{D})\in \mathcal{G}$ is a very ample linear series $\mathcal{D}$ on the curve $C=\xi^{-1}(p)$. We set $r:=h^0(C, |\mathcal{D}|)-1$ for a general $(p,\mathcal{D})\in{\mathcal{G}}$. By the Castelnuovo genus bound, we have $r\le 5$.
 
\vskip 4pt
\noindent
(4-a)
If $r=5$, the curve $C$ embedded by the complete linear series $|\mathcal{D}|$ in $\mathbb{P}^5$ is an extremal curve of degree $13$. There are three possibilities  for the residual series $\mathcal{F}=|K_C-\mathcal{D}|=g^3_9$; either $\mathcal{F}$  is (i) compounded, (ii) birationally very ample or (iii) very ample. If $\mathcal{F}$ is compounded, then 
$C$ is either trigonal or bi-elliptic with $\mathcal{F}$ having a base point. However if $C$ is bi-elliptic then the series $|K_C-\mathcal{F}|=|\mathcal{D}|$ cannot be very ample which one may verify rather easily.
Hence we have $\mathcal{F}=3g^1_3$. On the other hand, it is easy to check that  a trigonal curve $C$ of genus $g=12$  is embedded in $\mathbb{P}^5$ as a curve of degree $g+1$ by the series $|K_C-3g^1_3|$. Therefore it follows that $$\dim\mathcal{G}=\dim\mathcal{M}^1_{g,3}+\dim\mathbb{G}(4,5)=(2g+1)+5=30$$ and the component $\mathcal{H}_1$ arising in this way has dimension $\dim\mathcal{G}+\dim\mathbb{P}GL(5)=54$, which is the expected dimension. 
\vskip 4pt
\noindent
If $\mathcal{F}=g^3_9$ is birationally very ample or very ample, then $C$ can be embedded  in $\mathbb{P}^3$ as a curve smooth of type $(4,5)$ on a quadric. Hence the family consisting of such curves form a family $\mathcal{H}$ of dimension $\dim\mathbb{P}H^0(\mathbb{P}^3, \mathcal{O}(2))+\dim\mathbb{P}H^0(\mathbb{P}^1\times\mathbb{P}^1, \mathcal{O}(4,5))=9+29=38$ and therefore it follows that that
\begin{align*}
\dim\mathcal{G}&=\dim\mathcal{W}^\vee +\dim\mathbb{G}(4,5)=\dim\mathcal{H}-\dim\mathbb{P}GL(4)+\dim\mathbb{G}(4,5)\\&=28 <\lambda (13,12,4)=30,
\end{align*}
contrary to Propostion \ref{facts} (1);  here $\mathcal{W}^\vee\subset\mathcal{W}^3_{g-3}$ denotes the locus 
consisting  of the residual  series of  elements in $\mathcal{G}$. Hence the family of curves in $\mathbb{P}^4$ induced by a general $4$-dimensional subseries of a complete very ample series $|\mathcal{D}|=g^5_{13}$ whose residual series is birationally very ample or very ample does not constitute a component of $\mathcal{H}_{13,12,4}$. Therefore we conclude that the only possible component of $\mathcal{H}_{13,12,4}$ whose general element is not linearly normal is $\mathcal{H}_1$ consisting of trigonal curves.
\vskip 4pt
\noindent
(4-b)
If $r=4$, $|K_C-\mathcal{D}|=g^2_9$ is a net of degree $9$ for general $(p, \mathcal{D})\in \mathcal{G}\subset \mathcal{G}^{4}_{g+1}$. Let $\mathcal{W}^\vee\subset\mathcal{W}^2_{g-3}$ be the locus 
consisting  of the residual  series of  elements in $\mathcal{G}$. As in the previous case (4-a), there  are three possibilities for $|K_{C}-\mathcal{D}|=\mathcal{E}$; either $\mathcal{E}$  is (i) compounded, (ii) birationally very ample or (iii) very ample. $\mathcal{E}$ is not very ample since there is no integer solution for the equation $g=\frac{(e-1)(e-2)}{2}=12$ for any $e\le 9$. By applying Proposition \ref{wrdbd} to the locus $\mathcal{W}^\vee$ (for $d'=g-3$ and $r'=2$), together with the inequality $\dim\mathcal{G}=\dim\mathcal{W}^\vee \ge\lambda(13,12,4)$, we may exclude the possibility for $\mathcal{E}$ being compounded and conclude that 
$\mathcal{E}$ is birationally very ample,  base-point-free ($b=0$) and  $\dim\mathcal{G}=\lambda(13,12,4)$. Similar to the case $g=10$ in (3), we see that the component $\mathcal{G}$ is birational to the irreducible locus $\mathcal{G}'\subset\mathcal{G}^2_9$, which is the locus over which the Severi variety $\Sigma_{9,12}$ lies. Therefore we come up with another component $\mathcal{H}_2$ of the expected dimension 
$$\dim\mathcal{G}+\dim\mathbb{P}GL(5)=\lambda (13,12,4)+\dim\mathbb{P}GL(5)=54$$ such that a general element of $\mathcal{H}_2$ corresponds to a linearly normal curve and the residual series of the hyperplane series induces a plane model of degree $9$ with nodal singularities. 
Finally we remark that 
$\mathcal{H}_2\neq\emptyset$  since there exists a smooth curve in $\mathbb{P}^4$ of genus $12$ and degree $13$ which belongs to $\mathcal{H}_2$; 
cf. \cite[Theorem 1.0.1]{rathman}. Note  the either of $\mathcal{H}_i$'s is not contained in the other since both have the same dimension and hence $\mathcal{H}_1\neq\mathcal{H}_2$. We finally have $\mathcal{H}_{13,12,4}^\mathcal{L}=\mathcal{H}_2$ is irreducible whereas $\mathcal{H}_{13,12,4}=\mathcal{H}_1\cup\mathcal{H}_2$ is not.
\qed

\begin{rmk}\label{genus9} (i) The two components $\mathcal{H}_1$ and $\mathcal{H}_2$ of $\mathcal{H}_{13,12,4}$  have the same expected dimension, which is a quite rare case; one component with linearly normal curves and the other one non-linearly normal curves.  Another example of a Hilbert scheme with two components of the same expected dimension is $\mathcal{H}_{9,10,3}$. A general element of one component of $\mathcal{H}_{9,10,3}$ corresponds to a curve of type $(3,6)$ on a smooth quadric surface in $\mathbb{P}^3$ and the other component consist of curves which are complete intersection of two cubic surfaces. However a general element of any of these two components of $\mathcal{H}_{9,10,3}$ is a linearly normal curve unlike $\mathcal{H}_{13,12,4}$.

\vskip 4pt
\noindent
\rm(ii) The proof we presented for  Theorem \ref{lowgenus} (3) - the irreducibility of $\mathcal{H}_{g+1,g,4}$ when $g=10$ -  indicates the method of a proof we will follow for general cases.  
\end{rmk}

By Theorem \ref{lowgenus} we may assume that $g\ge 11$ and $g\neq 12$ for the rest of this section. 

\begin{lem}\label{dualbirva}
Let ${\mathcal{G}}\subset \widetilde{\mathcal{G}}_\mathcal{L}\subset \widetilde{\mathcal{G}}\subset\mathcal{G}^{4}_{g+1}$ be an irreducible component whose general element $(p, \mathcal{D})$ is a very ample and complete linear series $\mathcal{D}$ on the curve $C=\xi^{-1}(p)$. Then 
\begin{enumerate}
\item[\rm{(1)}] $\dim{\mathcal{G}}=4g-18$.
\item[\rm{(2)}] a general element of the component ${\mathcal{W}}^{\vee}\subset \mathcal{W}^2_{g-3}$ consisting of the residual  series of elements in  ${\mathcal{G}}$ is a base-point-free, complete  and birationally very ample net.

\end{enumerate}
\end{lem}
\noindent
{\it Proof.}
\quad By Proposition \ref{facts}(1), we have
\begin{equation}\label{min}
\lambda (g+1,g,4)=3g-3+\rho(g+1,g,4)=4g-18\le \dim {\mathcal{G}}.
\end{equation}
Note that  $\dim\mathcal{D}=\dim|\mathcal{D}|=4$ for a general $(p,\mathcal{D})\in{\mathcal{G}}$ and let $\mathcal{W}\subset \mathcal{W}^{4}_{g+1}$ be the component containing the image of the natural rational map 
${\mathcal{G}}\overset{\iota}{\dashrightarrow} \mathcal{W}^{4}_{g+1}$ with $\iota (\mathcal{D})=|\mathcal{D}|$.
We also let $\mathcal{W}^\vee\subset \mathcal{W}^{2}_{g-3}$ be the locus consisting  of the residual  series  of elements in $\mathcal{W}$, i.e. $\mathcal{W}^\vee =\{(p, \omega_C\otimes L^{-1}): (p, L)\in\mathcal{W}\}.$
\vskip 4pt
\noindent
(a)
If a general element of $\mathcal{W}^{\vee}$ is compounded, then by Proposition \ref{wrdbd}(c),
\begin{eqnarray*}
4g-18&\le&\dim\mathcal{G}\le\dim \mathcal{W}=\dim \mathcal{W}^{\vee}\\
&\le& 2g-1+(g-3)-2\cdot 2\\
&=&3g-8
\end{eqnarray*}
implying  $g\le 10$ contrary to our assumption $g\ge 11$. 
Therefore we conclude that a general element of $\mathcal{W}^{\vee}$ is either very ample or birationally very ample. 
\vskip 4pt
\noindent
(b)
Suppose that the moving part of a  general element of $\mathcal{W}^{\vee}\subset\mathcal{W}^{2}_{g-3}$ is very ample and let $b$ be the degree of the base locus $B$ of a general element of $\mathcal{W}^{\vee}$. Putting $e:=g-3-b$ and $g=\frac{(e-1)(e-2)}{2}\ge 11$,  a simple numerical calculation yields $e\ge 7$ and

\begin{align*}\dim\mathcal{G}&=\dim\mathbb{P}(H^0(\mathbb{P}^2, \mathcal{O}(e)))-\dim\mathbb{P}\rm{GL}(3)+b\\&=\frac{(e+1)(e+2)}{2}-9+b
<\lambda(g+1,g,4)=4g-18,
\end{align*}
contrary to the inequality (\ref{min}).
\vskip 4pt
\noindent
(c)
Therefore the moving part of a  general element of $\mathcal{W}^{\vee}\subset\mathcal{W}^{2}_{g-3}$ is birationally very ample and we let $b$ be the degree of the base locus $B$ of a general element of $\mathcal{W}^{\vee}$. 
 By Proposition \ref{wrdbd}(b), we have
\begin{eqnarray*}
4g-18&\le&\dim\mathcal{G}\le\dim \mathcal{W}= \dim \mathcal{W}^{\vee}\\
&\le& 3(g-3)+g-1-4\cdot 2-2b\\
&=&4g-18-2b, 
\end{eqnarray*}
implying $b=0$. Thus it finally follows that  

\begin{equation}\label{equal}
	\dim\mathcal{G}=\dim\mathcal{W}=\dim\mathcal{W}^{\vee}=4g-18
\end{equation}
as we wanted.
\qed
\vskip 4pt
The irreducibility of $\mathcal{H}^\mathcal{L}_{g+1,g,4}$ readily follows as an
immediate consequence of Lemma \ref{dualbirva} together with Lemma \ref{Gisirred}.
\vskip 4pt
\noindent
{\it Proof of Theorem \ref{H_g+1,g,4}.}
Retaining the same notations as before, let $\widetilde{\mathcal{G}}_\mathcal{L}$ be the union of irreducible component $\mathcal{G}$ of $\mathcal{G}^{4}_{g+1}$ whose general element corresponds to a pair $(p,\mathcal{D})$ such that $\mathcal{D}$ is very ample and complete  linear series on $C:=\xi^{-1}(p)$. Let $\widetilde{\mathcal{W}}^\vee$ be the union of the  components $\mathcal{W}^\vee$ of $\mathcal{W}^2_{g-3}$, where $\mathcal{W}^\vee$ consists of the residual series of elements in $\mathcal{G}\subset\tilde{\mathcal{G}}_\mathcal{L}$.
We also let $\mathcal{G}'$ be the union of irreducible components of $\mathcal{G}^{2}_{g-3}$ whose general element is a birationally very ample and  base-point-free linear series.  We recall that, by Lemma \ref{Gisirred} and Proposition \ref{facts}(2), $\mathcal{G'}$ is irreducible and $\dim \mathcal{G'}=3(g-3)+g-9=4g-18$.  
By Lemma \ref{dualbirva} (or by the equality (\ref{equal})), \begin{equation}\label{dominant}\dim \mathcal{W}^{\vee}=\dim\mathcal{G}=4g-18=\dim \mathcal{G}'.\end{equation}
Since a general element of any component $\mathcal{W}^{\vee}\subset\widetilde{\mathcal{W}}^\vee\subset\mathcal{W}^2_{g-3}$ is a base-point-free,  birationally very ample and complete net by Lemma \ref{dualbirva}, there is a natural rational map $\widetilde{\mathcal{W}}^\vee\overset{\kappa}{\dashrightarrow}\mathcal{G}'$ with $\kappa(|\mathcal{D}|)=\mathcal{D}$   which is clearly injective on an open subset  $\widetilde{\mathcal{W}}^{\vee o}$ of $\widetilde{\mathcal{W}}^\vee$ consisting of those which are base-point-free,  birationally very ample and complete nets. Therefore the  rational map $\kappa$ is dominant by (\ref{dominant}).
We also note that there is a natural rational map 
$\mathcal{G}' \overset{\iota}{\dashrightarrow}\widetilde{\mathcal{W}}{^\vee}$ with $\iota (\mathcal{D})=|\mathcal{D|}$, which is an inverse to $\kappa$ (wherever it is defined).
Therefore it follows that $\widetilde{\mathcal{W}}^{\vee}$ is birationally equivalent to  the irreducible locus $\mathcal{G}'$, hence $\widetilde{\mathcal{W}}^{\vee}$ is irreducible and so is $\widetilde{\mathcal{G}}_\mathcal{L}$. Since $\mathcal{H}_{g+1,g,4}^\mathcal{L}$ is a $\mathbb{P}\textrm{GL}(5)$-bundle over an open subset of $\widetilde{\mathcal{G}}_\mathcal{L}$,  $\mathcal{H}_{g+1,g,4}^\mathcal{L}$ is irreducible.\qed

\begin{rmk}
(i)
For $g\geq 9$, the Hilbert scheme $\mathcal{H}^\mathcal{L}_{g+1,g,4}$ of linearly normal curves is non-empty.
For low $g$ which has not been treated in  Theorem \ref{lowgenus}, e.g. $g=11$, one may argue as follows. 
Let $S$ be a del-Pezzo surface in $\mathbb{P}^4$ by blowing up $\mathbb{P}^2$ at $5$ points $p_1,\cdots , p_5$ in general position. 
The proper transformation of a plane septic with four nodal singularities at $p_1, \cdots ,p_4$ and passing through the remaining  point $p_5$ gives rise to smooth 
curve of genus $g=11$ and degree $12$ in $\mathbb{P}^4$, which is linearly normal.
In  the Brill-Noether range $g\ge 15$,   $\mathcal{H}^\mathcal{L}_{g+1,g,4}$ is non-empty which is  guaranteed by the existence of the principal component $\mathcal{G}_0$; cf. Remark \ref{principal} (2). 

\vskip 4pt
\noindent
(ii)
It is worthwhile to note that the Hilbert scheme $\mathcal{H}^\mathcal{L}_{g+1,g,4}$ of linearly normal curves is generically reduced. For $g=9$, the generically reducedness of $\mathcal{H}^\mathcal{L}_{10,9,4}$ follows from \cite[Corollory (3.1)]{CC}.
For $g\ge 10$, the generically reducedness of $\mathcal{H}^\mathcal{L}_{g+1,g,4}$ follows from the fact that the dimension of the singular locus of $\mathcal{G}'\subset \mathcal{G}^{2}_{g-3}$ does not exeed $g-8<\lambda(g+1,g,4)=\lambda(g-3,g,2)=4g-18$; cf. \cite[Proposition (2.9)]{AC2}. Thus, $\mathcal{G}'$ is generically reduced and since $\mathcal{G}_\mathcal{L}$ is birational to $\mathcal{G}'$, $\mathcal{H}_{g+1,g,4}^\mathcal{L}$,  is also generically reduced.

\vskip 4pt
\noindent
(iii) We expect that $\mathcal{H}_{g+1,g,4}=\mathcal{H}^\mathcal{L}_{g+1,g,4}$ except for $g=12$ so that 
$\mathcal{H}_{g+1,g,4}$ is irreducible unless $g=9$ or $g=12$. The issue is that one need to come up with an auxiliary result such that 
a general element of a component $\mathcal{G}$  of $\widetilde{\mathcal{G}}$ is complete and then a similar proof would work for the rest.
\end{rmk}

\bibliographystyle{spmpsci} 

\end{document}